\documentclass{article}

\usepackage[english,portuges]{babel}

\usepackage{amsmath,amsfonts,amssymb,amsthm}
\usepackage{graphicx}
\usepackage{url}

\usepackage{pstricks,pst-node,pst-text,pst-3d,pst-plot}
\usepackage{pst-circ}
\usepackage{pstricks-add}

\newtheorem{thm}{Proposi\c{c}\~{a}o}
\newtheorem{definition}[thm]{Defini\c{c}\~{a}o}
\newtheorem{ex}[thm]{Exemplo}
\newtheorem{obs}[thm]{Observa\c{c}\~{a}o}
\newtheorem{theorem}[thm]{Teorema}
\newtheorem{proposition}[thm]{Proposi\c{c}\~{a}o}

\newlength{\MiniPageLeft}
\newlength{\MiniPageRight}

\begin{document}

\selectlanguage{portuges}

\title{Escalas Temporais e \textsf{Mathematica}\thanks{Aceite
para publica\c{c}\~{a}o no Boletim da SPM
(Sociedade Portuguesa de Matem\'{a}tica).}}

\author{Artur M. C. Brito da Cruz$^1$\\
\url{acruz@est.ips.pt}
\and
Helena Sofia Rodrigues$^2$\\
\url{sofiarodrigues@esce.ipvc.pt}
\and
Delfim F. M. Torres$^3$\\
\url{delfim@ua.pt}}

\date{$^1$Escola Superior de Tecnologia\\
Instituto Polit\'{e}cnico de Set\'{u}bal\\[0.3cm]
$^2$Escola Superior de Ci\^{e}ncias Empresariais\\
Instituto Polit\'{e}cnico de Viana do Castelo\\[0.3cm]
$^3$Departamento de Matem\'{a}tica\\
Universidade de Aveiro}

\maketitle


\begin{abstract}
As escalas temporais s\~{a}o um modelo de tempo,
onde os casos cl\'{a}ssicos de tempo cont\'{\i}nuo e discreto
s\~{a}o considerados e fundidos num mesmo quadro geral.
Neste artigo apresentam-se as defini\c{c}\~{o}es b\'{a}sicas
em escalas temporais e introduz-se, paralelamente,
um \emph{package} em \textsf{Mathematica}.

\smallskip

\noindent \textbf{Palavras chave:} escalas temporais,
\textsf{Mathematica}.
\end{abstract}


\selectlanguage{english}

\begin{abstract}
Time scales are a model of time, where the continuous
and the discrete time cases are considered and merged
into the same framework. In this paper some basic definitions
of the time scale calculus are presented. Simultaneously,
a package in \textsf{Mathematica} is introduced.

\smallskip

\noindent \textbf{Keywords:} time scales,
computer algebra system \textsf{Mathematica}.

\smallskip

\noindent \textbf{2010 Mathematics Subject Classification:}
26-04, 26E70.
\end{abstract}

\selectlanguage{portuges}


\section{Introdu\c{c}\~{a}o}

De acordo com E. T. Bell (1883-1960) ``uma das tarefas principais da Matem\'{a}tica
\'{e} harmonizar o cont\'{\i}nuo e o discreto incluindo-os numa matem\'{a}tica abrangente
e eliminando a obscuridade de ambos'' \cite{Bell}.
Em 1988 Stefan Hilger introduziu o c\'{a}lculo em escalas temporais --- \emph{Time Scales} \cite{Hilger1988}
--- e iniciou-se uma rela\c{c}\~{a}o prof\'{\i}cua entre o c\'{a}lculo diferencial e o c\'{a}lculo
\`{a}s diferen\c{c}as \cite{Hilger1990}. O \emph{C\'{a}lculo em Escalas Temporais} introduzido por Hilger
permite cumprir com sucesso a tarefa proposta por Bell.
\'{E} um assunto que tem recebido nos \'{u}ltimos anos particular aten\c{c}\~{a}o
e onde se tem vindo a assistir a interessantes avan\c{c}os.

O formalismo das escalas temporais
tem um enorme potencial para aplica\c{c}\~{o}es
em diversas \'{a}reas tais como a biologia, a teoria do controlo, a economia e a medicina,
onde os sistemas din\^{a}micos envolvidos cont\^{e}m frequentemente
uma parte discreta e uma parte cont\'{\i}nua
\cite{Atici2008,Biles2005,Ferreira2008}. Por exemplo, um consumidor recebe
o sal\'{a}rio num dado momento do m\^{e}s (tempo discreto), mas vai ponderando ao longo de todo o tempo quanto
deve gastar e quanto deve poupar (tempo cont\'{\i}nuo) \cite{Atici2008,Biles2005}.
A biologia \'{e} fecunda em exemplos, pois \'{e} normal o crescimento de plantas
e insectos depender fortemente de uma \'{e}poca do ano, devido a factores como a temperatura, humidade
e pluviosidade \cite{Duke2006}.
Por exemplo, \'{e} natural que se considere um modelo din\^{a}mico
em escalas temporais para estudar insectos \emph{Magicicada},
que exibem uma combina\c{c}\~{a}o de ciclos de vida longos e curtos.
Na verdade, as cigarras passam v\'{a}rios anos de crescimento subterr\^{a}neo
como jovens (de 4 a 17 anos, dependendo da esp\'{e}cie),
saindo depois acima do solo apenas por um curto est\'{a}gio adulto
de v\'{a}rias semanas. Desta forma faz todo o sentido considerar-se
uma escala de tempo diferente para os per\'{\i}odos ``abaixo'' e ``acima do solo''.
Outro exemplo \'{e} dado pelo insecto \emph{Magicicada septendecim}
que vive como uma larva por 17 anos e como adulto por cerca de uma semana.
Na medicina, Jones \emph{et al.} \cite{Jones2004} apresentam
uma modela\c{c}\~{a}o em escalas temporais para usar o desbridamento de uma ferida como
um controlo na cicatriza\c{c}\~{a}o natural. Em Matem\'{a}tica encontram-se tamb\'{e}m v\'{a}rias aplica\c{c}\~{o}es,
nomeadamente na \'{a}rea das desigualdades matem\'{a}ticas \cite{GronTS,si:ts,jia,[1]SpecialIssue:Aveiro2},
na teoria do controlo \cite{Bartosiewicz2007,B:Paw,Bartosiewicz2008,withEwaP:avoidance,Gosia:Delfim},
c\'{a}lculo das varia\c{c}\~{o}es \cite{comRicardoISO:nabla,Atici2008,Bartos:Del,Ferreira2008,Basia:post_doc_Aveiro:1,NataliaHigherOrderNabla}
e optimiza\c{c}\~{a}o multi-objectivo \cite{specialAveiro2Basia}.

Na Sec\c{c}\~{a}o~\ref{sec:2} introduzem-se os conceitos b\'{a}sicos relativos
\`{a}s escalas temporais e mostra-se como esta teoria \'{e} uma generaliza\c{c}\~{a}o
do caso real (tempo cont\'{\i}nuo) e do caso dos n\'{u}meros inteiros (tempo discreto).
Na Sec\c{c}\~{a}o~\ref{sec:3} apresentam-se duas no\c{c}\~{o}es de diferenciabilidade em escalas temporais:
as derivadas delta e as derivadas nabla. O texto \'{e} acompanhado com exemplos de
utiliza\c{c}\~{a}o do sistema de computa\c{c}\~{a}o alg\'{e}brica \textsf{Mathematica}. Termina-se
com a Sec\c{c}\~{a}o~\ref{sec:4} de conclus\~{a}o e notas finais sobre algumas
linhas de investiga\c{c}\~{a}o nesta \'{a}rea.


\section{No\c{c}\~{o}es elementares}
\label{sec:2}

Uma \emph{escala temporal}, $\mathbb{T}$, \'{e} um subconjunto n\~{a}o vazio fechado de $\mathbb{R}$.
Os conjuntos
\begin{center}
\begin{tabular}{l}
$\mathbb{T}    \mathbb{=}\mathbb{R}$,\\
$\mathbb{T}    \mathbb{=}\mathbb{Z}$,\\
$\mathbb{T}    \mathbb{=}h\mathbb{Z} =\left\{  hk:k\in\mathbb{Z}\right\}  \text{ com }h>0$,\\
$\mathbb{T}    \mathbb{=}\overline{q^{\mathbb{Z}}}
=\left\{  q^{k}:k\in \mathbb{Z}\right\}  \cup\left\{  0\right\}  \text{ com }q>1$,\\
$\mathbb{T}    \mathbb{=}\mathbb{P}_{a,b}=\underset{k=0}{\overset{+\infty}{\cup}}\left[
k\left(  a+b\right)  ,k\left(  a+b\right)  +a\right]  \text{ onde }a,b>0$,
\end{tabular}
\end{center}
s\~{a}o exemplos de escalas temporais.

Neste trabalho usa-se o \emph{package} em \textsf{Mathematica}
\emph{TimeScales} que foi criado para o efeito em colabora\c{c}\~{a}o com
o Prof. Pedro A. F. Cruz (\url{pedrocruz@ua.pt}) da Universidade de Aveiro.
O leitor \'{e} convidado a efectuar o \emph{download} do \emph{package} em
\url{http://www.esce.ipvc.pt/docentes/srodrigues/TimeScales.m}
e a fazer as suas pr\'{o}prias experi\^{e}ncias. Para isso, \'{e} necess\'{a}rio colocar
o ficheiro na directoria \emph{Extrapackages} do \textsf{Mathematica} e,
na janela de trabalhos do \textsf{Mathematica}, executar o comando
\texttt{Needs["TimeScales`"]}. Para introduzir uma escala temporal em \textsf{Mathematica}
basta depois consider\'{a}-la como uma lista, onde se pode introduzir valores discretos e cont\'{\i}nuos.
Note-se que os valores devem ser escritos por ordem crescente e,
caso haja a uni\~{a}o de v\'{a}rios subconjuntos,
estes devem ser disjuntos dois a dois.


\bigskip
\bigskip

\begin{minipage}[c]{0.38\linewidth}

\bigskip

\noindent
\footnotesize{
$ts:=\{Dados\}$

onde os valores em ``Dados'' podem ser discretos,
em forma de intervalo, ou a uni\~{a}o dos dois casos anteriores.}
\end{minipage}\hspace*{\fill}
\begin{minipage}[c]{0.6\linewidth}
\centering
\includegraphics[scale=0.4]{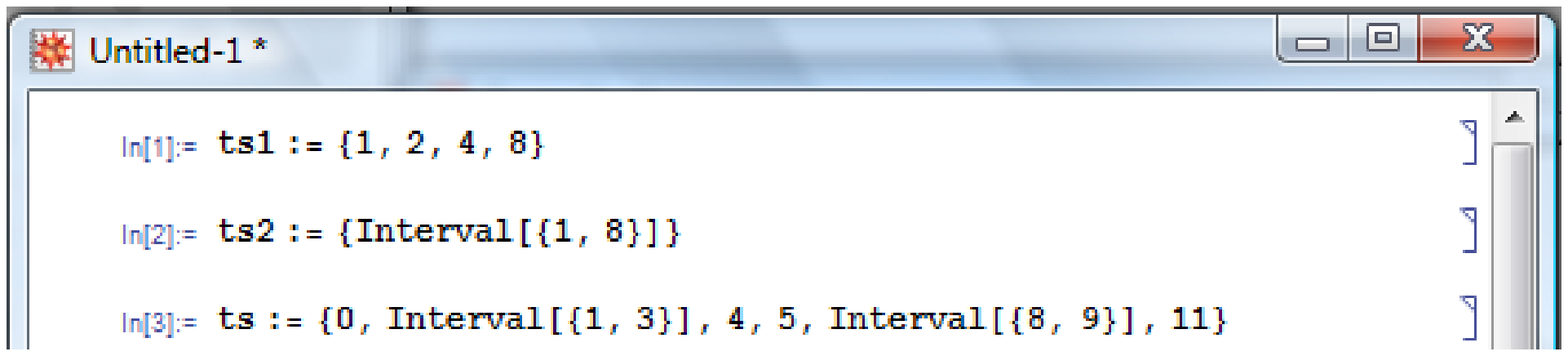}
\end{minipage}

\medskip
\bigskip
\bigskip


Neste artigo considerar-se-\'{a} a seguinte escala temporal como exemplo ilustrativo:
\texttt{ts:=\{0, Interval[\{1,3\}],4,5,Interval[\{8,9\}],11\}}, o que na nota\c{c}\~{a}o usual
da Matem\'{a}tica nada mais \'{e} que o conjunto $ts =\{ 0, 4, 5, 11 \} \cup [1, 3] \cup [8, 9]$.

Para uma visualiza\c{c}\~{a}o em \textsf{Mathematica} dos elementos que pertencem
\`{a} escala temporal usa-se o seguinte comando:


\smallskip
\smallskip

\begin{minipage}[c]{0.25\linewidth}
\noindent
\footnotesize{
$TSPlot[escalatemporal]$}
\end{minipage}\hspace*{\fill}
\begin{minipage}[c]{0.58\linewidth}
\centering
\includegraphics[scale=0.4]{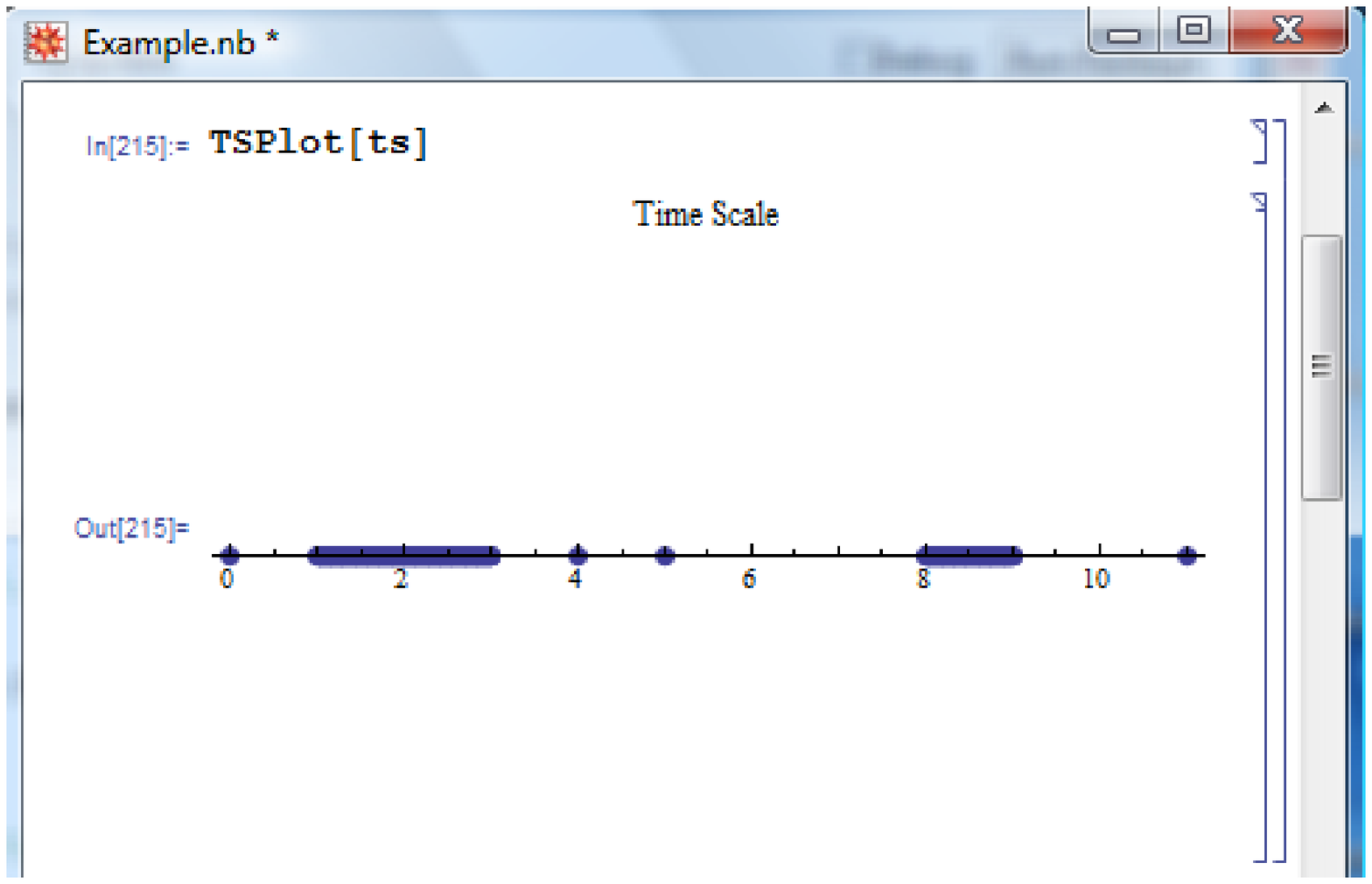}
\end{minipage}

\smallskip
\smallskip


Dada uma escala temporal $\mathbb{T}$ e $t\in\mathbb{T}$,
definem-se os seguintes operadores:

\begin{itemize}
\item o \textit{operador de avan\c{c}o }$\sigma:\mathbb{T\rightarrow T}$ por
$\sigma\left(  t\right)  =\inf\left\{  s\in\mathbb{T}:s>t\right\}$,
onde $\sigma\left(  \sup\mathbb{T}\right)  =\sup\mathbb{T}$ caso
$\sup\mathbb{T<+\infty}$;

\item o \textit{operador de recuo} $\rho:\mathbb{T\rightarrow T}$ por
$\rho\left(  t\right)  =\sup\left\{  s\in\mathbb{T}:s<t\right\}$,
onde $\rho\left(  \inf\mathbb{T}\right)  =\inf\mathbb{T}$ caso $\inf
\mathbb{T>-\infty}$;

\item a \textit{fun\c{c}\~{a}o de rarefac\c{c}\~{a}o de avan\c{c}o} $\mu:\mathbb{T\rightarrow
}\left[  0,+\infty\right[$ por
$\mu\left(  t\right)  =\sigma\left(  t\right)  -t$;

\item a \textit{fun\c{c}\~{a}o de rarefac\c{c}\~{a}o de recuo} $\nu:\mathbb{T\rightarrow
}\left[  0,+\infty\right[$ por
$\nu\left(  t\right)  =t-\rho\left(  t\right)$.
\end{itemize}

\begin{ex}
A Figura~\ref{fig:op:av:rec} ilustra o comportamento dos operadores de avan\c{c}o
e recuo numa escala de tempo $\mathbb{T}$.
\begin{figure}[!ht]
\centering
\[
\begin{pspicture}(0,-0.6)(\textwidth,1.8)

\rput[l](0.1,0.5){$\mathbb{T}$}
\psline[linewidth=1.5pt]{-}(0.08\textwidth,0.5)(0.25\textwidth,0.5)
\pscircle*(0.2\textwidth,0.5){0.05}
\pscircle*(0.3\textwidth,0.5){0.05}
\pscircle*(0.4\textwidth,0.5){0.05}
\pscircle*(0.45\textwidth,0.5){0.05}
\psline[linewidth=1.5pt]{-}(0.5\textwidth,0.5)(0.72\textwidth,0.5)
\pscircle*(0.72\textwidth,0.5){0.05}
\pscircle*(0.75\textwidth,0.5){0.05}
\pscircle*(0.8\textwidth,0.5){0.05}
\pscircle*(0.83\textwidth,0.5){0.05}
\pscircle*(0.86\textwidth,0.5){0.05}
\pscircle*(0.91\textwidth,0.5){0.05}
\psline[linewidth=1.5pt]{-}(0.91\textwidth,0.5)(\textwidth,0.5)

\psline[linewidth=0.2pt]{->}(0.2\textwidth,-0.2)(0.2\textwidth,0.47)
\rput[c](0.2\textwidth,-0.4){{\small $\rho(t_0)=t_0=\sigma(t_0)$}}

\psline[linewidth=0.2pt]{->}(0.3\textwidth,1.2)(0.3\textwidth,0.53)
\rput(0.3\textwidth,1.4){{\small $\rho(t_1)$}}

\psline[linewidth=0.2pt]{->}(0.4\textwidth,1.2)(0.4\textwidth,0.53)
\rput(0.4\textwidth,1.4){{\small $t_1$}}

\psline[linewidth=0.2pt]{->}(0.45\textwidth,-0.2)(0.45\textwidth,0.47)
\rput(0.45\textwidth,-0.4){{\small $\sigma(t_1)$}}


\psline[linewidth=0.2pt]{->}(0.72\textwidth,1.2)(0.72\textwidth,0.53)
\rput[rr](0.72\textwidth,1.4){{\small $\rho(t_2)=t_2$}}

\psline[linewidth=0.2pt]{->}(0.75\textwidth,-0.2)(0.75\textwidth,0.47)
\rput(0.75\textwidth,-0.4){{\small $\sigma(t_2)$}}


\psline[linewidth=0.2pt]{->}(0.86\textwidth,1.2)(0.86\textwidth,0.53)
\rput(0.86\textwidth,1.4){{\small $\rho(\!t_3\!)$}}

\psline[linewidth=0.2pt]{->}(0.91\textwidth,1.2)(0.91\textwidth,0.53)
\rput[l](0.9\textwidth,1.4){{\small $t_3\!\!=\!\!\sigma\!(\!t_3\!)$}}

\end{pspicture}
\]
\caption{Ilustra\c{c}\~{a}o dos operadores de avan\c{c}o e recuo.}
\label{fig:op:av:rec}
\end{figure}
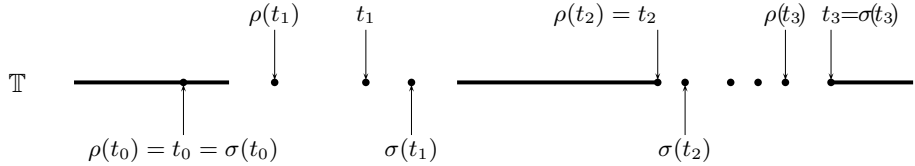
\end{ex}

Em \textsf{Mathematica}:


\bigskip
\bigskip

\begin{minipage}[c]{0.42\linewidth}
\noindent
\footnotesize{
$TSsigma[escalatemporal, ponto]$

caso se pretenda calcular o operador de avan\c{c}o em apenas um ponto;

\medskip

$\{\#, TSsigma[escalatemporal, \#]\} \& / \text{@}$

$Range[Min[escala temporal], $

$Max[escalatemporal], passo]$

caso se pretenda calcular o operador de avan\c{c}o numa s\'{e}rie de pontos
(note-se que caso os pontos n\~{a}o perten\c{c}am \`{a} escala temporal,
o \textsf{Mathematica} emite um aviso);

\medskip

$TSsigmaPlot[escalatemporal]$

representa\c{c}\~{a}o gr\'{a}fica do operador de avan\c{c}o.}

\hspace*{0.1cm}
\end{minipage}
\begin{minipage}[c]{0.58\linewidth}
\centering
\includegraphics[scale=0.4]{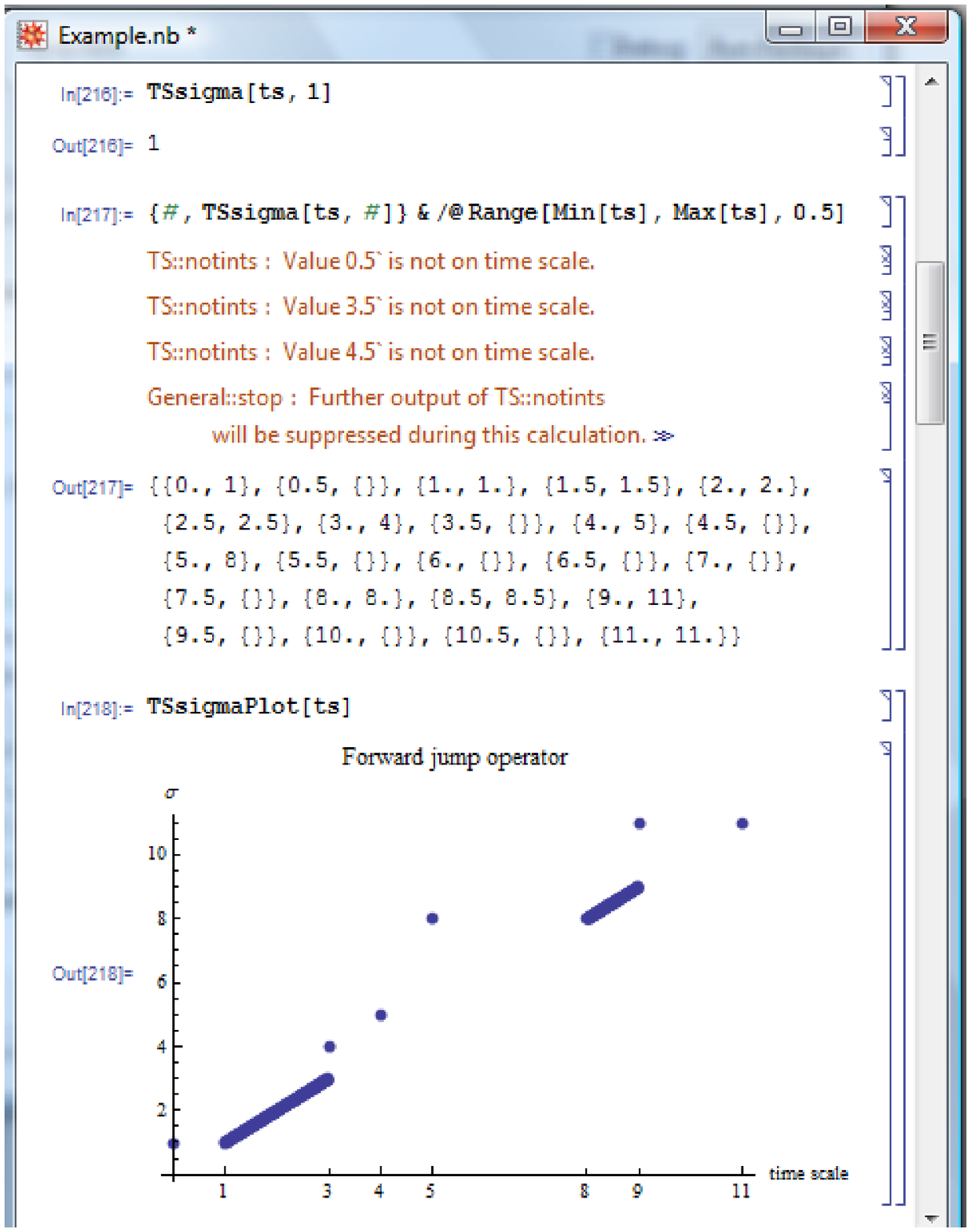}
\end{minipage}

\bigskip


Os comandos \textsf{Mathematica} relativos ao operador de recuo e \`{a}s fun\c{c}\~{o}es de rarefac\c{c}\~{a}o s\~{a}o,
respectivamente,
\begin{verbatim}
         TSrho[escalatemporal, ponto]
         TSmu[escalatemporal, ponto]
         TSnu[escalatemporal, ponto]
\end{verbatim}

Os operadores de avan\c{c}o e recuo permitem classificar os pontos de uma escala
temporal. Um ponto $t\in\mathbb{T}$ diz-se:

\begin{itemize}
\item \textit{discreto \`{a} direita} se $\sigma\left(  t\right)  >t$;

\item \textit{denso \`{a} direita} se $\sigma\left(  t\right)  =t$;

\item \textit{denso \`{a} esquerda} se $\rho\left(  t\right)  =t$;

\item \textit{discreto \`{a} esquerda} se $\rho\left(  t\right)  <t$.
\end{itemize}

Um ponto denso \`{a} direita e denso \`{a} esquerda diz-se simplesmente \textit{denso};
um ponto diz-se \textit{isolado} caso seja discreto \`{a} esquerda e \`{a} direita simultaneamente.

\begin{ex}
Na Figura~\ref{fig:op:av:rec} tem-se que
$t_{0}$ \'{e} denso; $t_{1}$ \'{e} isolado;
$t_{2}$ \'{e} denso \`{a} esquerda e discreto \`{a} direita;
e $t_{3}$ \'{e} discreto \`{a} esquerda e denso \`{a} direita.
\end{ex}

\begin{ex}
Se $\mathbb{T}=\mathbb{R}$, ent\~{a}o todo o ponto $t\in\mathbb{T}$ \'{e} denso pois
$\sigma\left(  t\right)  =t=\rho\left(  t\right)$. Al\'{e}m disso,
$\mu\left(t\right) \equiv 0$ e $\nu\left(  t\right) \equiv 0$.
\end{ex}

\begin{ex}
Se $\mathbb{T}=\mathbb{Z}$, ent\~{a}o todo o ponto $t\in\mathbb{T}$ \'{e} isolado
pois $\sigma\left(  t\right)=t+1$ e $\rho\left(  t\right)  =t-1$. Al\'{e}m disso,
$\mu\left(  t\right) \equiv 1$ e $\nu\left(  t\right)  \equiv 1$.
\end{ex}

\begin{ex}
Se $\mathbb{T=}h \mathbb{Z}$, $h>0$, ent\~{a}o todo o ponto $t\in\mathbb{T}$ \'{e} isolado pois
\begin{align*}
\sigma\left(  t\right)   &= \inf\left\{  s\in\mathbb{T}:s>t\right\}
= t+h\\
\rho\left(  t\right)   &= \sup\left\{  s\in\mathbb{T}:s<t\right\}
=t-h.
\end{align*}
Neste caso $\mu\left(  t\right)  =h$ e $\nu\left(  t\right)  =h$.
\end{ex}

\begin{ex}
Se $\mathbb{T=}\overline{q^{\mathbb{Z}}}$, $q > 1$,
ent\~{a}o para $t=q^{n_{0}}\in\mathbb{T}$ tem-se
\begin{align*}
\sigma\left(  t\right)   &= \inf\left\{  s\in\mathbb{T}:s>t\right\}
=q^{n_{0}+1}=qt\\
\rho\left(  t\right)   &= \sup\left\{  s\in\mathbb{T}:s<t\right\}
=q^{n_{0}-1}=q^{-1}t
\end{align*}
e, para $t=0$,
\[
\sigma\left(  0\right)  =0=\rho\left(  0\right) \, .
\]
Logo todo o ponto $t\neq0$ \'{e} isolado e o ponto $t=0$ \'{e} denso. A fun\c{c}\~{a}o de
rarefac\c{c}\~{a}o de avan\c{c}o \'{e} dada por
\[
\mu\left(  t\right)  =\left\{
\begin{array}[c]{cll}
\left(  q-1\right)  t & \text{ se } & t\in\mathbb{T}\backslash\left\{
0\right\} \\
0 & \text{ se } & t=0 \, ;
\end{array}
\right.
\]
a fun\c{c}\~{a}o de rarefac\c{c}\~{a}o de recuo \'{e} dada por
\[
\nu\left(  t\right)  =\left\{
\begin{array}[c]{cll}
\left(  1-q^{-1}\right)  t & \text{ se } & t\in\mathbb{T}\backslash\left\{
0\right\} \\
0 & \text{ se } & t=0 \, .
\end{array}
\right.
\]
\end{ex}

Nos exemplos anteriores $\sigma$ e $\rho$ s\~{a}o fun\c{c}\~{o}es inversas:
$\sigma = \rho^{-1}$ e $\rho = \sigma^{-1}$. Em geral,
conforme se ilustra no Exemplo~\ref{ex7},
$\sigma \ne \rho^{-1}$ e $\rho \ne \sigma^{-1}$.

\begin{ex}
\label{ex7}
Sejam $a,b>0$ e $\mathbb{T=P}_{a,b}=
\underset{k=0}{\overset{+\infty}{\cup}}\left[  k\left(  a+b\right), k\left(  a+b\right)  +a\right]$.
Ent\~{a}o,
\[
\sigma\left(  t\right)  =\left\{
\begin{array}[c]{cll}
t & \text{ se } & t\in\underset{k=0}{\overset{+\infty}{\cup}}\left[  k\left(
a+b\right), k\left(  a+b\right)  +a\right[ \\
t+b & \text{ se } & t\in\underset{k=0}{\overset{+\infty}{\cup}}\left\{
k\left(  a+b\right)  +a\right\}
\end{array}
\right.
\]
\[
\unitlength 1mm 
\linethickness{0.2pt}
\ifx\plotpoint\undefined\newsavebox{\plotpoint}\fi 
\begin{picture}(13,25)(0,-13)
\put(-15,-15){\vector(0,2){30}}

\linethickness{1pt}
\put(-15,-15){\line(1,0){6}}
\put(-6,-15){\line(1,0){6}}
\put(3,-15){\line(1,0){6}}
\put(12,-15){\line(1,0){6}}
\put(16,-15){\vector(1,0){2}}
\put(-15,-17){\makebox(0,0)[cc]{\scriptsize $0$}}
\put(-9,-17){\makebox(0,0)[cc]{\scriptsize $a$}}
\put(17,-17){\makebox(0,0)[cc]{\scriptsize $t$}}
\put(-15.5,-9){\line(1,0){1}}
\put(-15.5,-6){\line(1,0){1}}
\put(-15.5,0){\line(1,0){1}}
\put(-15.5,3){\line(1,0){1}}
\put(-15.5,9){\line(1,0){1}}
\put(-15.5,12){\line(1,0){1}}
\put(-17,-9){\makebox(0,0)[rc]{\scriptsize $a$}}
\put(-17,-6){\makebox(0,0)[rc]{\scriptsize $a+b$}}
\put(-17,0){\makebox(0,0)[rc]{\scriptsize $2a+b$}}

{\black
\thicklines
\put(-15,-15){\line(1,1){5.55}}
\put(-6,-6){\line(1,1){5.55}}
\put(3,3){\line(1,1){5.55}}
\put(12,12){\line(1,1){4}}
\put(-9,-9){\circle{1}}
\put(-9,-6){\circle*{1}}
\put(-6,-6){\circle*{1}}
\put(0,0){\circle{1}}
\put(0,3){\circle*{1}}
\put(3,3){\circle*{1}}
\put(9,9){\circle{1}}
\put(9,12){\circle*{1}}
\put(12,12){\circle*{1}}
\put(12,3){\makebox(0,0)[cc]{$\sigma(t)$}}
}
\end{picture}
\]

\[
\rho\left(  t\right)  =\left\{
\begin{array}
[c]{cll}%
t & \text{~se~} & t\in\underset{k=0}{\overset{+\infty}{\cup}}\left]  k\left(
a+b\right)  ,k\left(  a+b\right)  +a\right]  \cup\left\{  0\right\} \\
t-b & \text{~se~} & t\in\underset{k=1}{\overset{+\infty}{\cup}}\left\{
k\left(  a+b\right)  \right\}
\end{array}
\right.
\]

\[
\unitlength 1mm 
\linethickness{0.2pt}
\ifx\plotpoint\undefined\newsavebox{\plotpoint}\fi 
\begin{picture}(13,31)(0,-17)
\put(-15,-15){\vector(0,2){30}}

\linethickness{0.5pt}
\put(-15,-15){\line(1,0){6}}
\put(-6,-15){\line(1,0){6}}
\put(3,-15){\line(1,0){6}}
\put(12,-15){\line(1,0){6}}
\put(16,-15){\vector(1,0){2}}
\put(-15,-17){\makebox(0,0)[cc]{\scriptsize $0$}}
\put(-9,-17){\makebox(0,0)[cc]{\scriptsize $a$}}
\put(17,-17){\makebox(0,0)[cc]{\scriptsize $t$}}
\put(-15.5,-9){\line(1,0){1}}
\put(-15.5,-6){\line(1,0){1}}
\put(-15.5,0){\line(1,0){1}}
\put(-15.5,3){\line(1,0){1}}
\put(-15.5,9){\line(1,0){1}}
\put(-15.5,12){\line(1,0){1}}
\put(-17,-9){\makebox(0,0)[rc]{\scriptsize $a$}}
\put(-17,-6){\makebox(0,0)[rc]{\scriptsize $a+b$}}
\put(-17,0){\makebox(0,0)[rc]{\scriptsize $2a+b$}}

{\black
\thicklines
\put(-15,-15){\line(1,1){5.6}}
\put(-5.55,-5.55){\line(1,1){5.6}}
\put(3.45,3.45){\line(1,1){5.6}}
\put(12.45,12.45){\line(1,1){4}}
\put(-9,-9){\circle*{1}}
\put(-6,-9){\circle*{1}}
\put(-6,-6){\circle{1}}
\put(0,0){\circle*{1}}
\put(3,0){\circle*{1}}
\put(3,3){\circle{1}}
\put(9,9){\circle*{1}}
\put(12,9){\circle*{1}}
\put(12,12){\circle{1}}
\put(12,3){\makebox(0,0)[cc]{$\rho(t)$}}
}

\end{picture}
\]

\[
\mu\left(  t\right)  =\left\{
\begin{array}[c]{cll}
0 & \text{ se } & t\in\underset{k=0}{\overset{+\infty}{\cup}}\left[  k\left(
a+b\right)  ,k\left(  a+b\right)  +a\right[ \\
b & \text{ se } & t\in\underset{k=0}{\overset{+\infty}{\cup}}\left\{  k\left(
a+b\right)  +a\right\}
\end{array}
\right.
\]

\[
\unitlength 1mm 
\linethickness{0.2pt}
\ifx\plotpoint\undefined\newsavebox{\plotpoint}\fi 
\begin{picture}(13,15)(0,-10)
\put(-15,-6){\vector(0,1){13}}

\linethickness{0.5pt}
\put(-15,-6){\line(1,0){5.55}}
\put(-6,-6){\line(1,0){5.55}}
\put(3,-6){\line(1,0){5.55}}
\put(12,-6){\line(1,0){5.55}}
\put(16,-6){\vector(1,0){2}}
\put(-15,-8){\makebox(0,0)[cc]{\scriptsize $0$}}
\put(-9,-8){\makebox(0,0)[cc]{\scriptsize $a$}}
\put(17,-8){\makebox(0,0)[cc]{\scriptsize $t$}}
\put(-15.5,-3){\line(1,0){1}}
\put(-17,-3){\makebox(0,0)[rc]{\scriptsize $b$}}

{\thicklines
\black
\put(-15,-6){\line(1,0){5.55}}
\put(-9,-6){\circle{1}}
\put(-9,-3){\circle*{1}}
\put(-6,-6){\line(1,0){5.55}}
\put(0,-6){\circle{1}}
\put(-6,-6){\circle*{1}}
\put(3,-6){\line(1,0){5.55}}
\put(0,-3){\circle*{1}}
\put(3,-6){\circle*{1}}
\put(12,-6){\line(1,0){5.55}}
\put(9,-6){\circle{1}}
\put(9,-3){\circle*{1}}
\put(12,-6){\circle*{1}}
}

{\black
\put(9,3){\makebox(0,0)[cc]{$\mu(t)$}}
}
\end{picture}
\]
\[
\nu\left(  t\right)  =\left\{
\begin{array}[c]{cll}
0 & \text{ se } & t\in\underset{k=0}{\overset{+\infty}{\cup}}\left]  k\left(
a+b\right), k\left(  a+b\right)  +a\right]  \cup\left\{  0\right\} \\
b & \text{ se } & t\in\underset{k=1}{\overset{+\infty}{\cup}}\left\{  k\left(
a+b\right)  \right\}
\end{array}
\right.
\]
Para $k\geqslant0$, visto que
\[
\rho\left(  k\left(  a+b\right)  +a\right)  =k\left(  a+b\right)  +a
\text{ \ e \ }\sigma\left(  k\left(  a+b\right)  +a\right)  =\left(  k+1\right)
\left(  a+b\right) ,
\]
ent\~{a}o $\sigma\circ\rho\neq id$,
onde $id$ designa a aplica\c{c}\~{a}o identidade. Para $k\geqslant 1$ tem-se que
\[
\sigma\left(  k\left(  a+b\right)  \right)  =k\left(  a+b\right)  \text{ \ e
\ }\rho\left(  k\left(  a+b\right)  \right)  =k\left(  a+b\right)  -b ,
\]
logo $\rho\circ\sigma\neq id$.
\end{ex}

\begin{proposition}
Um ponto $t\in\mathbb{T}$ \'{e} discreto \`{a} direita e denso \`{a} esquerda se e s\'{o} se
$\sigma\circ\rho\left(  t\right)  \neq t$.
\end{proposition}

\begin{proof}
Repare-se que se $t$ \'{e} um ponto discreto \`{a} esquerda, ent\~{a}o $\rho\left(
t\right)  =t_{1}\neq t$ e $\sigma\left(  t_{1}\right)  =t$ e, portanto,
$\sigma\circ\rho\left(  t\right)  =t$. Assim a desigualdade $\sigma\circ
\rho\left(  t\right)  \neq t$ garante que $t$ \'{e} denso \`{a} esquerda e, como tal,
\[
\sigma\circ\rho\left(  t\right)  =\sigma\left(  t\right)  \neq t \, ,
\]
ou seja, $t$ \'{e} um ponto denso \`{a} esquerda e discreto \`{a} direita.
A implica\c{c}\~{a}o contr\'{a}ria \'{e} imediata.
\end{proof}

\begin{obs}
\begin{enumerate}
\item Um ponto $t\in\mathbb{T}$ \'{e} denso \`{a} direita e discreto \`{a} esquerda se e
s\'{o} se $\rho\circ\sigma\left(  t\right)  \neq t$.

\item O operador de avan\c{c}o \'{e} sobrejectivo se e s\'{o} se n\~{a}o existem pontos em
$\mathbb{T}$ que sejam simultaneamente discretos \`{a} direita e densos \`{a} esquerda.

\item O operador de avan\c{c}o \'{e} injectivo se e s\'{o} se n\~{a}o existem pontos em
$\mathbb{T}$ que sejam simultaneamente densos \`{a} direita e discretos \`{a} esquerda.
\end{enumerate}
\end{obs}

Para se proceder \`{a} classifica\c{c}\~{a}o de um ponto numa escala temporal
em \textsf{Mathematica} basta introduzir o comando


\bigskip
\bigskip

\begin{minipage}[c]{0.40\linewidth}
\noindent
\footnotesize{
$TSPoint[escalatemporal, ponto]$

}

\end{minipage}\hspace*{\fill}
\begin{minipage}[c]{0.58\linewidth}
\centering
\includegraphics[scale=0.4]{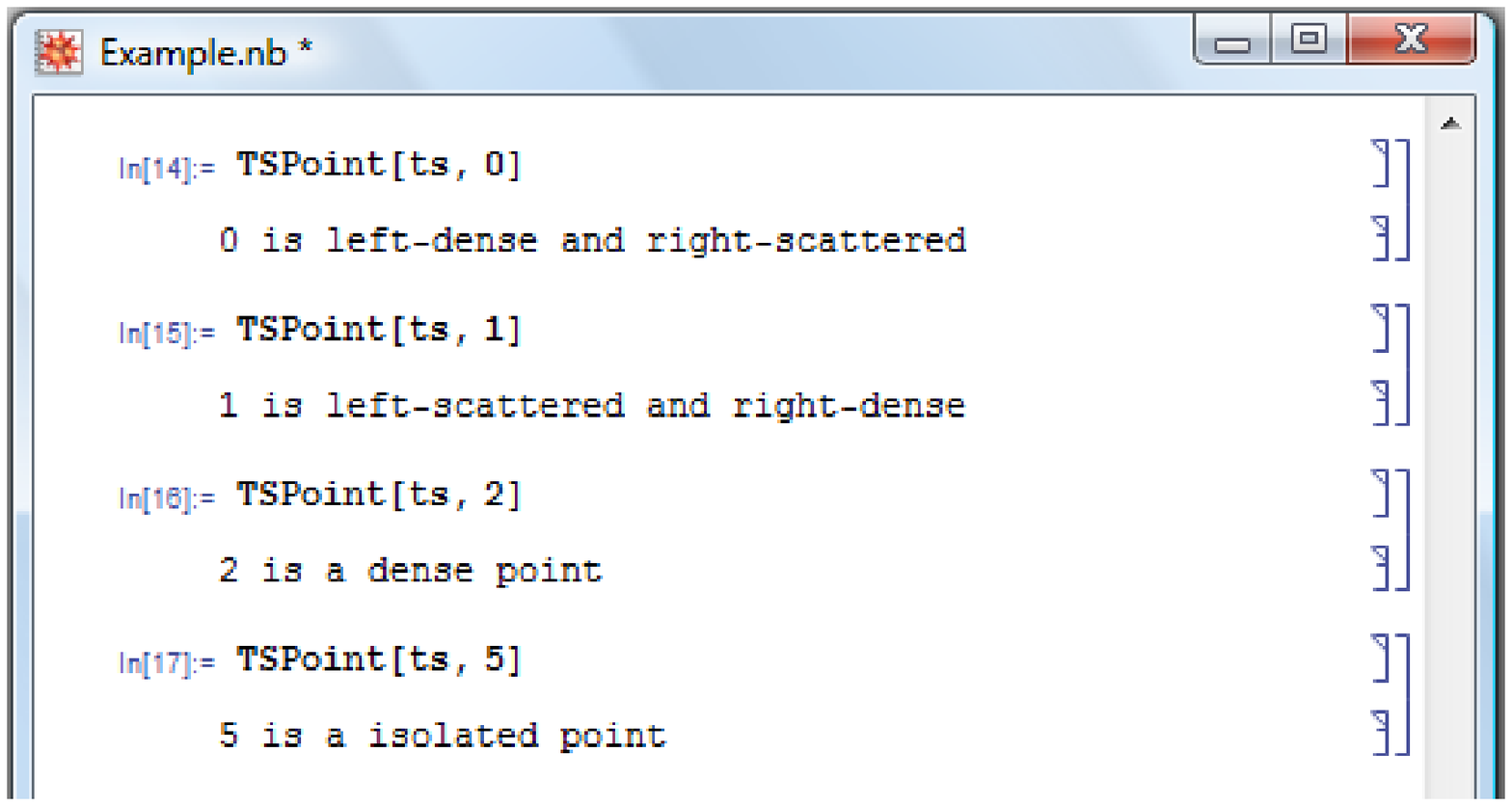}
\end{minipage}


\section{Diferenciabilidade}
\label{sec:3}

Seja $\mathbb{T}$ uma escala temporal e considere-se em $\mathbb{T}$ a
topologia induzida pela topologia usual dos n\'{u}meros reais.
De modo semelhante ao tempo discreto $\mathbb{T}=\mathbb{Z}$,
onde \'{e} usual considerar-se dois operadores de diferen\c{c}as finitas,
\begin{equation*}
\begin{split}
\Delta f(t_i) &= f(t_{i+1})-f(t_i) \, , \\
\nabla f(t_i) &= f(t_i)-f(t_{i-1}) \, ,
\end{split}
\end{equation*}
tamb\'{e}m na teoria das escalas temporais s\~{a}o comuns
duas no\c{c}\~{o}es de diferencia\c{c}\~{a}o (e integra\c{c}\~{a}o \cite{Malinowska2009}).


\subsection{Derivada delta}

De modo a introduzir-se a defini\c{c}\~{a}o de derivada delta \'{e} necess\'{a}rio considerar um
novo conjunto, $\mathbb{T}^{\kappa}$, definido do seguinte modo:
\[
\mathbb{T}^{\kappa}=\left\{
\begin{array}
[c]{cll}
\mathbb{T}\backslash\left\{  \sup\mathbb{T}\right\}  & \text{ } & \text{ se
}\rho\left(  \sup\mathbb{T}\right)  <\sup\mathbb{T<+\infty},\\
\mathbb{T} & \text{ } & \text{ caso contr\'{a}rio.}
\end{array}
\right.
\]

\begin{definition}
\label{def:der:delta}
Seja $f:\mathbb{T\rightarrow} \mathbb{R}$ uma fun\c{c}\~{a}o e seja $t\in\mathbb{T}^{\kappa}$. A \textit{derivada delta} de $f$ em
$t$, que se representa por $f^{\Delta}\left(  t\right)$, \'{e} o n\'{u}mero real
(caso exista) tal que para qualquer $\varepsilon>0$ existe uma vizinhan\c{c}a $U_{t}^{\delta}$
de $t$ em $\mathbb{T}$ com
\begin{equation}
\label{eq:def:der:delta}
\left\vert \left[  f\left(  \sigma\left(  t\right)  \right)  -f\left(
s\right)  \right]  -f^{\Delta}\left(  t\right)  \left[  \sigma\left(
t\right)  -s\right]  \right\vert \leqslant\varepsilon\left\vert \sigma\left(
t\right)  -s\right\vert
\end{equation}
para todo o $s\in U_{t}^{\delta}$.
Diz-se que $f$ \'{e} delta-diferenci\'{a}vel em $\mathbb{T}$ se
existe derivada delta de $f$ em $t$ para todo o
$t\in\mathbb{T}^{\kappa}$. A fun\c{c}\~{a}o $f^{\Delta}:\mathbb{T}^{\kappa}\mathbb{\rightarrow} \mathbb{R}$
diz-se a derivada delta de $f$ em $\mathbb{T}^{\kappa}$.
\end{definition}

\begin{theorem}
\label{thm:unic:delta}
Seja $f : \mathbb{T}\rightarrow \mathbb{R}$ e
$t \in \mathbb{T}^\kappa$. Ent\~{a}o $f$ tem no m\'{a}ximo uma derivada delta em $t$.
\end{theorem}
\begin{proof}
Suponha-se que $f$ tem duas derivadas $x, y \in \mathbb{R}$ em $t$.
Por defini\c{c}\~{a}o de $\mathbb{T}^\kappa$, toda a vizinhan\c{c}a de $t$
cont\'{e}m algum $s \in \mathbb{T}^\kappa$ com $s \ne \sigma(t)$.
Tem-se ent\~{a}o, para todo o $\varepsilon > 0$, um $s \ne \sigma(t)$ com
\begin{equation*}
\begin{split}
\left\vert   f\left(\sigma\left(t\right)  \right)
-f\left(s\right)  - x \cdot \left(  \sigma\left(
t\right)  -s\right)  \right\vert &\leqslant
\frac{\varepsilon}{2}\left\vert \sigma\left(t\right)  -s\right\vert \\
\left\vert f\left(\sigma\left(t\right)  \right)
-f\left(s\right)  - y \cdot \left(  \sigma\left(
t\right)  -s\right)  \right\vert &\leqslant
\frac{\varepsilon}{2}\left\vert \sigma\left(t\right)  -s\right\vert \, .
\end{split}
\end{equation*}
Isto implica que
\begin{equation*}
\begin{split}
|x - y| &\cdot |\sigma(t)-s| = | (x - y) \cdot (\sigma(t)-s) |\\
&\le \left| x \cdot (\sigma(t)-s) - \left[f(\sigma(t))-f(s)\right]\right|
+ \left| f(\sigma(t))-f(s) - y \cdot (\sigma(t)-s) \right|\\
&\le \varepsilon \cdot |\sigma(t)-s| \, .
\end{split}
\end{equation*}
Como $|\sigma(t)-s| \ne 0$, quando $\varepsilon$ tende para zero
conclui-se que $x = y$.
\end{proof}

\begin{obs}
Se $t \in \mathbb{T} \setminus \mathbb{T}^\kappa$, ent\~{a}o $f^\Delta(t)$
n\~{a}o est\'{a} definida de modo \'{u}nico. Com efeito, para tal $t$ pequenas
vizinhan\c{c}as $U_{t}^{\delta}$ de $t$ consistem apenas de $t$ e, al\'{e}m disso,
tem-se $\sigma(t) = t$. Por conseguinte a desigualdade
\eqref{eq:def:der:delta} \'{e} verificada para qualquer n\'{u}mero $f^\Delta(t)$ arbitr\'{a}rio.
\'{E} esta a raz\~{a}o porque na Defini\c{c}\~{a}o~\ref{def:der:delta} se exclui um ponto
maximal que seja discreto \`{a} esquerda.
\end{obs}

\begin{ex}
Seja $f:\mathbb{T\rightarrow} \mathbb{R}$, $\alpha$ uma constante
e $f\left(  t\right)  =\alpha$ para todo o
$t\in\mathbb{T}$. Ent\~{a}o $f^{\Delta}\left(  t\right) =0$
para todo o $t\in\mathbb{T}^{\kappa}$.
\end{ex}

\begin{ex}
Se $f:\mathbb{T\rightarrow}\mathbb{R}$ \'{e} a fun\c{c}\~{a}o identidade
$f\left(  t\right)  =t$, ent\~{a}o $f^{\Delta}\left(  t\right) =1$.
\end{ex}

\begin{ex}
Seja $f:\mathbb{T\rightarrow}
\mathbb{R}$ definida por $f\left(  t\right)  =t^{2}$.
Ent\~{a}o $f^{\Delta}\left(t\right) =t+\sigma\left(  t\right)=2t+\mu(t)$.
Para $\mathbb{T}=\mathbb{R}$ tem-se $f^{\Delta}=2t = f'(t)$;
para $\mathbb{T}=\mathbb{Z}$, $f^{\Delta}=2t+1 = \Delta f(t)$.
\end{ex}

\begin{obs}
Uma fun\c{c}\~{a}o $f:\mathbb{T\rightarrow} \mathbb{R}$ diz-se cont\'{\i}nua em $t_{0}$
se para qualquer $\varepsilon>0$ existir uma
vizinhan\c{c}a $U_{t_{0}}^{\delta} = \left] t_{0}-\delta,t_{0}+\delta\right[ \cap \mathbb{T}$
de $t_{0}$, $\delta>0$, tal que
$\left\vert f\left(  t\right)  -f\left(  t_{0}\right)  \right\vert <\varepsilon$
qualquer que seja o $t\in U_{t_{0}}^{\delta}$.
A fun\c{c}\~{a}o $f$ diz-se cont\'{\i}nua se for cont\'{\i}nua para todo o $t\in\mathbb{T}$.
Todas as fun\c{c}\~{o}es $f : \mathbb{T} \rightarrow \mathbb{R}$ s\~{a}o cont\'{\i}nuas
para $\mathbb{T} = \mathbb{Z}$; para $\mathbb{T}=\mathbb{R}$ tem-se
a defini\c{c}\~{a}o usual de continuidade em $\mathbb{R}$;
na escala temporal $\mathbb{T=P}_{1,1}$ as restri\c{c}\~{o}es das fun\c{c}\~{o}es reais
de vari\'{a}vel real cont\'{\i}nuas em $\left[  2k,2k+1\right]$,
$k\in\mathbb{N}_{0}$, s\~{a}o cont\'{\i}nuas, no entanto o operador de avan\c{c}o $\sigma$ n\~{a}o \'{e} cont\'{\i}nuo,
pois n\~{a}o \'{e} cont\'{\i}nuo nos pontos $2k+1$, $k\in \mathbb{N}_{0}$
(\textrm{cf.} Exemplo~\ref{ex7}).
\end{obs}

\begin{theorem}[\cite{Bohner2001}]
Seja $f:\mathbb{T\rightarrow}\mathbb{R}$
uma fun\c{c}\~{a}o e seja $t\in\mathbb{T}^{\kappa}$.

\begin{enumerate}
\item Se $f$ \'{e} delta-diferenci\'{a}vel em $t$, ent\~{a}o $f$ \'{e} cont\'{\i}nua em $t$.

\item Se $f$ \'{e} cont\'{\i}nua em $t$, com $t$ um ponto discreto \`{a} direita,
ent\~{a}o $f$ \'{e} delta-diferenci\'{a}vel em $t$ e
\[
f^{\Delta}\left(  t\right)  =\frac{f\left(  \sigma\left(  t\right)  \right)
-f\left(  t\right)  }{\mu\left(  t\right)  }\text{.}
\]

\item Se $t$ \'{e} denso \`{a} direita, ent\~{a}o $f$ \'{e} delta-diferenci\'{a}vel em $t$
se e s\'{o} se o limite
\[
\lim_{s\rightarrow t}\frac{f\left(  s\right)  -f\left(  t\right)  }{s-t}
\]
existe (e \'{e} finito). Nesse caso
\[
f^{\Delta}\left(  t\right)  =\lim_{s\rightarrow t}\frac{f\left(  s\right)
-f\left(  t\right)  }{s-t}\, .
\]

\item Se $f$ \'{e} delta-diferenci\'{a}vel em $t$, ent\~{a}o
\[
f\left(  \sigma\left(  t\right)  \right)  =f\left(  t\right)  +\mu\left(
t\right)  f^{\Delta}\left(  t\right) \, .
\]
\end{enumerate}
\end{theorem}

\begin{ex}
Se $\mathbb{T=R}$,
ent\~{a}o $f: \mathbb{R}\rightarrow \mathbb{R}$
\'{e} diferenci\'{a}vel no sentido delta em
$t\in \mathbb{R}$ se e s\'{o} se o limite
\[
\lim_{s\rightarrow t}\frac{f\left(  s\right)  -f\left(  t\right)  }{s-t}
\]
existe (e \'{e} finito), ou seja, se $f$ \'{e} diferenci\'{a}vel no sentido usual:
\[
f^{\Delta}\left(  t\right)  =\lim_{s\rightarrow t}\frac{f\left(  s\right)
-f\left(  t\right)  }{s-t}=f^{\prime}\left(  t\right)  \text{.}
\]
\end{ex}

\begin{ex}
Se $\mathbb{T}=\mathbb{Z}$, ent\~{a}o $f:\mathbb{Z}
\rightarrow \mathbb{R}$
\'{e} delta-diferenci\'{a}vel para todo o $t\in \mathbb{Z}$ e
\[
f^{\Delta}\left(  t\right)  =\frac{f\left(  \sigma\left(  t\right)  \right)
-f\left(  t\right)  }{\mu\left(  t\right)  }=f\left(  t+1\right)  -f\left(
t\right)  = \Delta f\left(  t\right)\, ,
\]
onde $\Delta$ \'{e} o operador usual de diferen\c{c}as finitas.
\end{ex}

\begin{ex}
Seja $\mathbb{T}= h \mathbb{Z}$, $h>0$. Uma vez que cada ponto $t\in\mathbb{T}$ \'{e} isolado,
\[
f^{\Delta}\left(  t\right)  =\frac{f\left(  \sigma\left(  t\right)  \right)
-f\left(  t\right)  }{\mu\left(  t\right)  }=\frac{f\left(  t+h\right)
-f\left(  t\right)  }{h} = \Delta_h f(t) \, .
\]
\end{ex}

\begin{ex}
A derivada usual do c\'{a}lculo-$q$ (qu\^{a}ntico) \cite{Kac},
tamb\'{e}m conhecida como derivada de Jackson \cite{Jackson},
\'{e} facilmente obtida escolhendo-se, para $q>1$,
$\mathbb{T}= q^{\mathbb{N}_0} := \{q^k : k \in \mathbb{N}_0\}$:
$$
f^\Delta(t) = \frac{f(q t) - f(t)}{(q-1) t} = \mathcal{D}_q f(t) \, .
$$
\end{ex}

\begin{ex}
Seja $\mathbb{T}\mathbb{=P}_{a,b}$ onde $a,b>0$. Considere-se uma
sucess\~{a}o de pontos $\left(t_{n}^{k}\right)_{n\in \mathbb{N}}
\subseteq \left[ k\left(  a+b\right)  ,k\left(  a+b\right)  +a\right[$
tal que
\[
\lim_{n\rightarrow+\infty}t_{n}^{k}=k\left(  a+b\right)  +a \, .
\]
Note-se que
\[
\lim_{n\rightarrow+\infty}\sigma\left(  t_{n}^{k}\right)
= \lim_{n\rightarrow +\infty}t_{n}^{k}
=k\left(  a+b\right)  +a\neq\sigma\left(  k\left(
a+b\right)  +a\right)  \, .
\]
Como o operador de avan\c{c}o n\~{a}o \'{e} cont\'{\i}nuo no ponto $t=k\left(  a+b\right)  +a$,
tem-se que $\sigma$ n\~{a}o \'{e} delta-diferenci\'{a}vel nesse ponto.
\end{ex}

Este \'{u}ltimo exemplo \'{e} um caso particular do resultado que se segue.

\begin{proposition}[\cite{Gosia:Delfim}]
Se $t\in\mathbb{T}^{\kappa}$, $t\neq\min\mathbb{T}$ e $\rho\left(  t\right)
=t<\sigma\left(  t\right)$, ent\~{a}o o operador de avan\c{c}o $\sigma$ n\~{a}o
\'{e} delta-diferenci\'{a}vel em $t$.
\end{proposition}

\begin{proof}
Com vista a um absurdo, considere-se que $\sigma$ \'{e} delta-diferenci\'{a}vel
em $t$ e que $\sigma^{\Delta}\left(  t\right) = a$. Ent\~{a}o para todo o $s\in
U_{t}$, vizinhan\c{c}a de $t$,
\[
\left\vert \left[  \sigma\left(  \sigma\left(  t\right)  \right)
-\sigma\left(  s\right)  \right]  -a\left[  \sigma\left(  t\right)  -s\right]
\right\vert \leqslant\varepsilon\left\vert \sigma\left(  t\right)
-s\right\vert \, .
\]
Em particular, para $s=t$
\[
\left\vert \left[  \sigma\left(  \sigma\left(  t\right)  \right)
-\sigma\left(  t\right)  \right]  -a\left[  \sigma\left(  t\right)  -t\right]
\right\vert \leqslant\varepsilon\left\vert \sigma\left(  t\right)
-t\right\vert \, .
\]
Ao tomar-se o limite quando $\varepsilon\rightarrow0$, obt\'{e}m-se:
\begin{equation*}
\left[  \sigma\left(  \sigma\left(  t\right)  \right)  -\sigma\left(
t\right)  \right]  -a\left[  \sigma\left(  t\right)  -t\right]  =0
 \Rightarrow a=\frac{\sigma\left(  \sigma\left(  t\right)  \right)
-\sigma\left(  t\right)  }{\sigma\left(  t\right)  -t} \, .
\end{equation*}
Por outro lado, visto que $t$ \'{e} discreto \`{a} direita e denso \`{a} esquerda, o ponto
$s\in U_{t}$ pode ser escolhido \`{a} esquerda de $t$. Logo, quando
$s\rightarrow t$, obt\'{e}m-se que $\sigma\left(  s\right) =s \rightarrow t$.
Ent\~{a}o,
\begin{align*}
& \left\vert \left[  \sigma\left(  \sigma\left(  t\right)  \right)  -t\right]
-a\left[  \sigma\left(  t\right)  -t\right]  \right\vert \leqslant
\varepsilon\left\vert \sigma\left(  t\right)  -t\right\vert \\
& \Rightarrow\left[  \sigma\left(  \sigma\left(  t\right)  \right)  -t\right]
-a\left[  \sigma\left(  t\right)  -t\right]  =0\\
& \Rightarrow a=\frac{\sigma\left(  \sigma\left(  t\right)  \right)
-t}{\sigma\left(  t\right)  -t} \, .
\end{align*}
Ao comparar os resultados obtidos para $a$, conclui-se que $\sigma\left(
t\right)  =t$, o que contradiz a hip\'{o}tese.
\end{proof}

Em \textsf{Mathematica} usa-se o comando \texttt{TSDelta}:


\medskip
\bigskip

\begin{minipage}[c]{0.45\linewidth}
\noindent
\footnotesize{
$TSDelta[escalatemporal, \text{\emph{fun\c{c}\~{a}o}}, ponto]$

Note-se que a fun\c{c}\~{a}o pode ser definida dentro ou fora do comando TSDelta.
}

\end{minipage}
\begin{minipage}[c]{0.5\linewidth}
\centering
\includegraphics[scale=0.45]{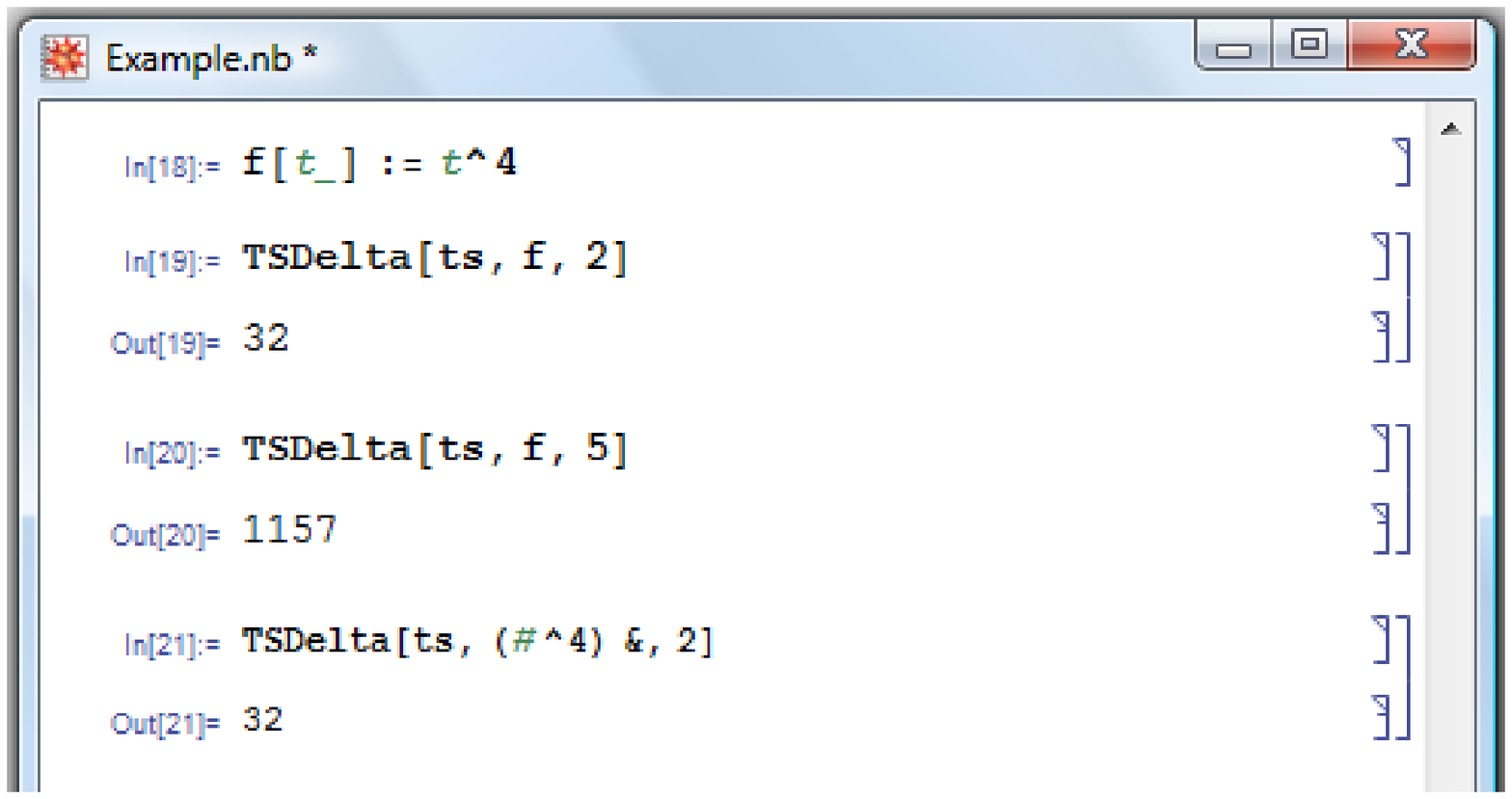}
\end{minipage}

\medskip
\bigskip


O teorema seguinte apresenta as propriedades fundamentais da derivada delta.

\begin{theorem}[\cite{Bohner2001}]
\label{thm:prp:fund:derDelta}
Sejam $f,g:\mathbb{T\rightarrow} \mathbb{R}$ duas fun\c{c}\~{o}es delta-diferenci\'{a}veis
em $t\in\mathbb{T}^{\kappa}$, $\alpha\in \mathbb{R}$. Ent\~{a}o,

\begin{enumerate}
\item $\left(  f+g\right)  ^{\Delta}\left(  t\right)  =f^{\Delta}\left(
t\right)  +g^{\Delta}\left(  t\right)$;

\item $\left(  \alpha f\right)  ^{\Delta}\left(  t\right)  =\alpha f^{\Delta
}\left(  t\right)$;

\item $\left(  fg\right)  ^{\Delta}\left(  t\right)  =f^{\Delta}\left(
t\right)  g\left(  t\right)  +f\left(  \sigma\left(  t\right)  \right)
g^{\Delta}\left(  t\right)  =f\left(  t\right)  g^{\Delta}\left(  t\right)
+f^{\Delta}\left(  t\right)  g\left(  \sigma\left(  t\right)  \right)$;

\item $\displaystyle\left(  \frac{f}{g}\right)  ^{\Delta}\left(  t\right)
=\frac{f^{\Delta}\left(  t\right)  g\left(  t\right)  -f\left(  t\right)
g^{\Delta}\left(  t\right)  }{g\left(  t\right)  g\left(  \sigma\left(
t\right)  \right)}$ se $g\left(  t\right)  g\left(  \sigma\left(  t\right)  \right)  \neq0$.
\end{enumerate}
\end{theorem}

\begin{obs}
Use-se a nota\c{c}\~{a}o $f^{\sigma}=f\circ\sigma$.
\begin{enumerate}
\item Sejam $f_{i}:\mathbb{T}\rightarrow \mathbb{R}$,
$i=1,\ldots, n$, fun\c{c}\~{o}es delta-diferenci\'{a}veis. Ent\~{a}o,
\begin{multline*}
\left(f_{1}f_{2}\cdots f_{n-1}f_{n}\right)^{\Delta}\\
=f_{1}^{\Delta}
f_{2} \cdots f_{n-1}f_{n}+f_{1}^{\sigma}f_{2}^{\Delta}f_{3} \cdots f_{n-1} f_{n}
+ \cdots +f_{1}^{\sigma}f_{2}^{\sigma} \cdots f_{n-1}^{\sigma}f_{n}^{\Delta}\, .
\end{multline*}

\item Seja $f:\mathbb{T}\rightarrow \mathbb{R}$ uma fun\c{c}\~{a}o delta-diferenci\'{a}vel. Ent\~{a}o,
\[
\left(  f^{n}\right)  ^{\Delta}=f^{\Delta}\left[  \sum_{k=1}^{n}f^{n-k}\left(
f^{\sigma}\right)^{k-1}\right]  \, .
\]

\item Seja $f:\mathbb{T}\rightarrow \mathbb{R}$
uma fun\c{c}\~{a}o delta-diferenci\'{a}vel em $t$ tal que $f\left(  t\right)  f\left(
\sigma\left(  t\right)  \right)  \neq 0$. Ent\~{a}o,
\[
\left(  f^{-n}\right)  ^{\Delta}
=-\frac{\left(  f^{n}\right)  ^{\Delta}}{f^{n}\left(  f^{\sigma}\right)^{n}}\, .
\]
\end{enumerate}
\end{obs}

As opera\c{c}\~{o}es do Teorema~\ref{thm:prp:fund:derDelta}
tamb\'{e}m est\~{a}o definidas no \emph{package}
\emph{TimeScales} do \textsf{Mathematica}:


\medskip
\bigskip

\begin{minipage}[c]{0.43\linewidth}
\noindent
\footnotesize{
$TSDeltaSum[escalatemporal, \text{\emph{fun\c{c}\~{a}o1}},$

$\text{\emph{fun\c{c}\~{a}o2}},ponto]$

\medskip

$TSDeltaProdConst[escalatemporal,$

$\text{\emph{fun\c{c}\~{a}o}},constante, ponto]$

\medskip

$TSDeltaProduct[escalatemporal,$

$\text{\emph{fun\c{c}\~{a}o1}}, \text{\emph{fun\c{c}\~{a}o2}},ponto]$

\medskip

$TSDeltaQuocient[escalatemporal,$

$\text{\emph{fun\c{c}\~{a}o1}}, \text{\emph{fun\c{c}\~{a}o2}},ponto]$

\medskip
}

\end{minipage}\hspace*{\fill}
\begin{minipage}[c]{0.58\linewidth}
\centering
\includegraphics[scale=0.45]{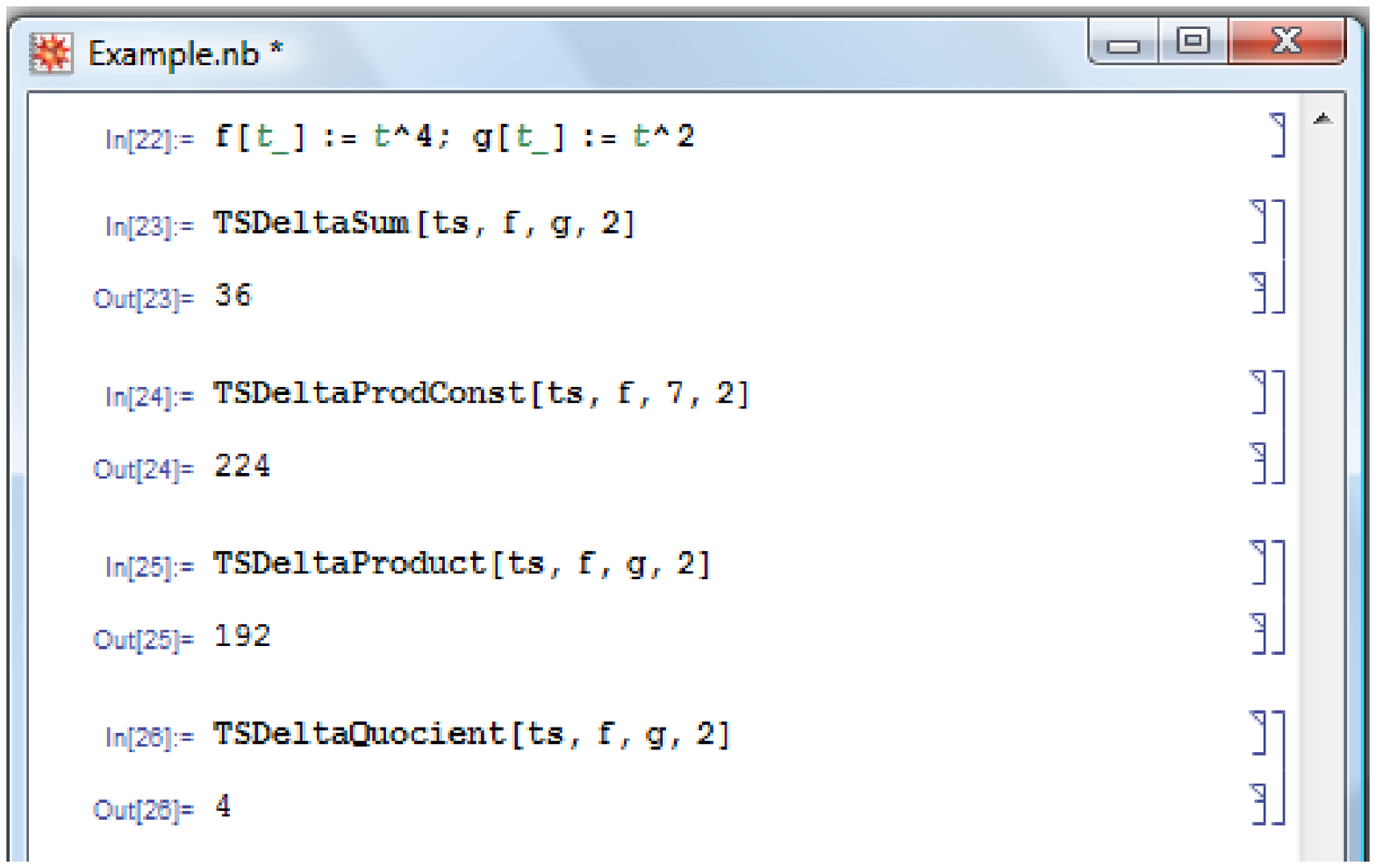}
\end{minipage}


\subsection{Derivada nabla}

Na sec\c{c}\~{a}o anterior definiu-se a derivada delta com recurso ao operador
de avan\c{c}o. Analogamente, introduz-se agora a no\c{c}\~{a}o de derivada nabla
recorrendo ao operador de recuo. Considere-se o seguinte conjunto:
\[
\mathbb{T}_{\kappa}=\left\{
\begin{array}
[c]{cll}
\mathbb{T}\backslash\left\{  \inf\mathbb{T}\right\}  & \text{} & \text{ se
}\mathbb{-\infty<}\inf\mathbb{T}<\sigma\left(  \inf\mathbb{T}\right), \\
\mathbb{T} & \text{} & \text{ caso contr\'{a}rio.}
\end{array}
\right.
\]

\begin{definition}
Seja $f:\mathbb{T}\rightarrow\mathbb{R}$ uma fun\c{c}\~{a}o e seja
$t\in\mathbb{T}_{\kappa}$. A \textit{derivada nabla} de $f$ em
$t$, que se representa por $f^{\nabla}\left(  t\right)  $, \'{e} o
n\'{u}mero real (caso exista) tal que para qualquer $\varepsilon>0$ existe uma
vizinhan\c{c}a $U_{t}^{\delta}$ de $t$ em $\mathbb{T}$ com
\[
\left\vert \left[  f\left(  \rho\left(  t\right)  \right)  -f\left(  s\right)
\right]  -f^{\nabla}\left(  t\right)  \left[  \rho\left(  t\right)
-s\right]  \right\vert \leqslant\varepsilon\left\vert \rho\left(  t\right)
-s\right\vert
\]
para todo o $s\in U_{t}^{\delta}$. Diz-se que $f$ \'{e} nabla-diferenci\'{a}vel
em $\mathbb{T}$ se existe derivada nabla de $f$ em $t$ para todo o
$t\in\mathbb{T}_{\kappa}$. A fun\c{c}\~{a}o $f^{\nabla}:\mathbb{T}_{\kappa}
\rightarrow \mathbb{R}$
diz-se a derivada nabla de $f$ em $\mathbb{T}_{\kappa}$.
\end{definition}

De modo semelhante ao caso delta (\textrm{cf.} Teorema~\ref{thm:unic:delta})
pode-se mostrar que se $f : \mathbb{T}\rightarrow \mathbb{R}$ e $t \in \mathbb{T}_\kappa$,
ent\~{a}o $f$ tem no m\'{a}ximo uma derivada nabla em $t$.

\begin{ex}
Para $\mathbb{T} = \mathbb{R}$ a derivada nabla coincide
com a no\c{c}\~{a}o usual de derivada: $f^{\nabla}=f^{\Delta}=f^{\prime}$.
\end{ex}

\begin{ex}
Para $\mathbb{T}=\mathbb{Z}$ tem-se
$f^{\nabla}\left(  t\right)  =\nabla f\left(  t\right)
=f\left(  t\right)  -f\left(  t-1\right)$.
\end{ex}

Os resultados seguintes s\~{a}o similares aos da derivada delta:

\begin{theorem}[\cite{Bohner2003}]
Seja $f:\mathbb{T} \rightarrow \mathbb{R}$
uma fun\c{c}\~{a}o e seja $t\in\mathbb{T}_{\kappa}$.

\begin{enumerate}
\item Se $f$ \'{e} nabla-diferenci\'{a}vel em $t$, ent\~{a}o $f$ \'{e} cont\'{\i}nua em $t$.

\item Se $f$ cont\'{\i}nua em $t$, com $t$ um ponto discreto \`{a} esquerda,
ent\~{a}o $f$ \'{e} nabla-diferenci\'{a}vel em $t$ e
\[
f^{\nabla}\left(  t\right)  =\frac{f\left(  t\right)  -f\left(
\rho\left(  t\right)  \right)  }{\nu\left(  t\right)  }\, .
\]

\item Se $t$ \'{e} denso \`{a} esquerda, ent\~{a}o $f$ \'{e} diferenci\'{a}vel
no sentido nabla em $t$ se e s\'{o} se o limite
\[
\lim_{s\rightarrow t}\frac{f\left(  s\right)  -f\left(  t\right)  }{s-t}
\]
existe (e \'{e} finito). Nesse caso,
\[
f^{\nabla}\left(  t\right)  =\lim_{s\rightarrow t}\frac{f\left(
s\right)  -f\left(  t\right)  }{s-t}\, .
\]

\item Se $f$ \'{e} nabla-diferenci\'{a}vel em $t$, ent\~{a}o
$f\left(  \rho\left(  t\right)  \right)  =f\left(  t\right)  -\nu\left(
t\right)  f^{\nabla}\left(  t\right)$.
\end{enumerate}
\end{theorem}

\begin{theorem}[\cite{Bohner2003}]
Sejam $f,g:\mathbb{T}\rightarrow \mathbb{R}$ duas
fun\c{c}\~{o}es nabla-diferenci\'{a}veis em $t\in\mathbb{T}_{\kappa}$,
$\alpha\in \mathbb{R}$. Ent\~{a}o,
\begin{enumerate}
\item $\left(  f+g\right)  ^{\nabla}\left(  t\right)
=f^{\nabla}\left(  t\right)  +g^{\nabla}\left(  t\right)$;

\item $\left(  \alpha f\right)  ^{\nabla}\left(  t\right)  =\alpha
f^{\nabla}\left(  t\right)$;

\item $\left(  fg\right)  ^{\nabla}\left(  t\right)
=f^{\nabla}\left(  t\right)  g\left(  t\right)  +f\left(
\rho\left(  t\right)  \right)  g^{\nabla}\left(  t\right)  =f\left(
t\right)  g^{\nabla}\left(  t\right)  +f^{\nabla}\left(
t\right)  g\left(  \rho\left(  t\right)  \right)  $;

\item $\displaystyle\left(  \frac{f}{g}\right)  ^{\nabla}\left(
t\right)  =\frac{f^{\nabla}\left(  t\right)  g\left(  t\right)
-f\left(  t\right)  g^{\nabla}\left(  t\right)  }{g\left(  t\right)
g\left(  \rho\left(  t\right)  \right)  }$
se $g\left(  t\right)  g\left(  \rho\left(  t\right)  \right)  \neq0$.
\end{enumerate}
\end{theorem}

Em \textsf{Mathematica} os comandos s\~{a}o an\'{a}logos aos da derivada delta
mas t\^{e}m sempre o prefixo \emph{TSNabla}.


\begin{minipage}[c]{0.45\linewidth}
\noindent
\footnotesize{$TSNabla[escalatemporal, \text{\emph{fun\c{c}\~{a}o}}, ponto]$

\smallskip

$TSNablaSum[escalatemporal, \text{\emph{fun\c{c}\~{a}o1}},$

$\text{\emph{fun\c{c}\~{a}o2}},ponto]$

\smallskip

$TSNablaProdConst[escalatemporal,$

$\text{\emph{fun\c{c}\~{a}o}}, constante, ponto]$

\smallskip

$TSNablaProduct[escalatemporal,$

$\text{\emph{fun\c{c}\~{a}o1}}, \text{\emph{fun\c{c}\~{a}o2}},ponto]$

\smallskip

$TSNablaQuocient[escalatemporal, $

$\text{\emph{fun\c{c}\~{a}o1}}, \text{\emph{fun\c{c}\~{a}o2}},ponto]$}

\end{minipage}\hspace*{\fill}
\begin{minipage}[c]{0.55\linewidth}
\centering
\includegraphics[scale=0.38]{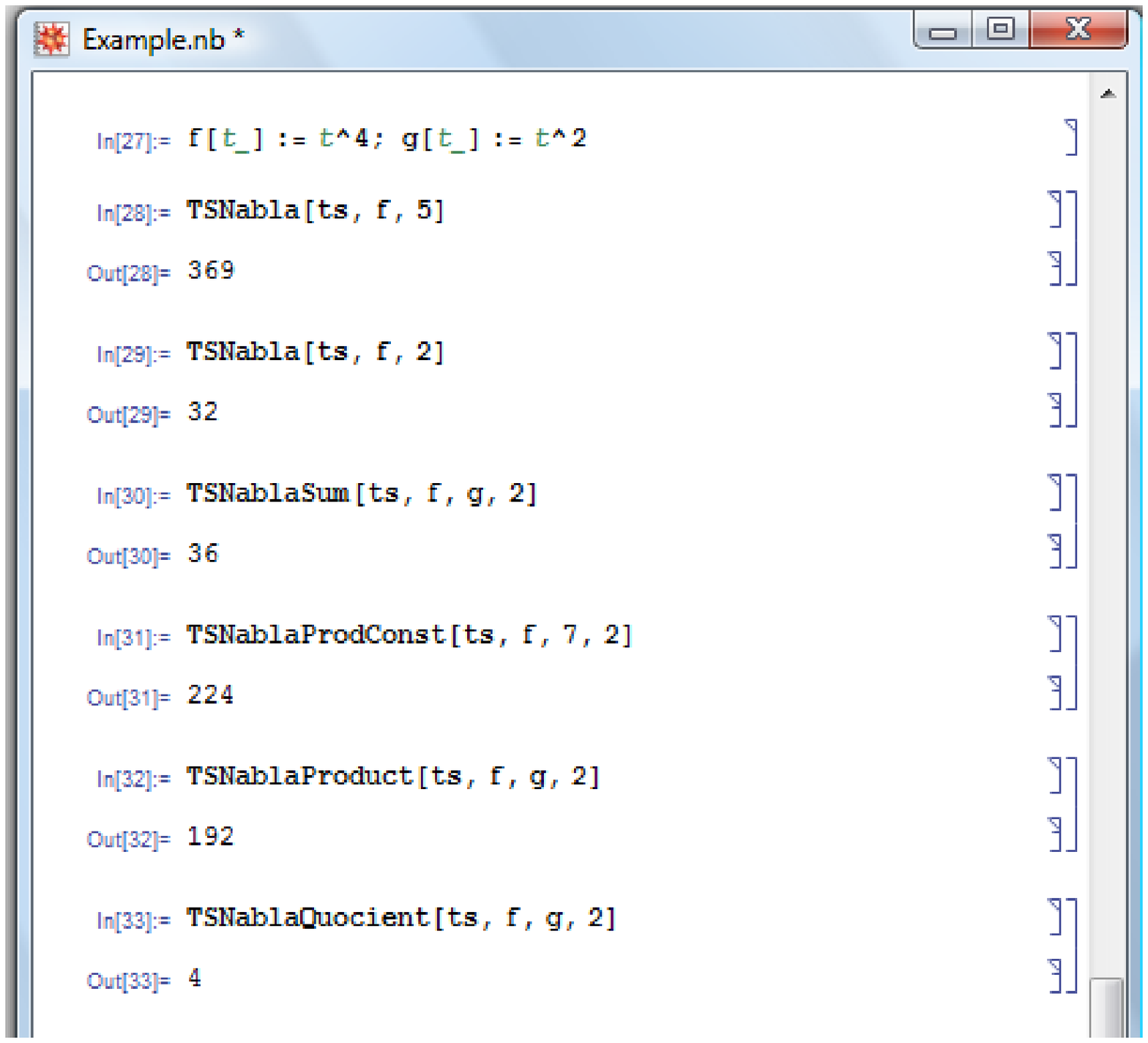}
\end{minipage}


\subsection{Regras da cadeia}

A f\'{o}rmula usual da derivada da fun\c{c}\~{a}o composta n\~{a}o \'{e} v\'{a}lida em todas as
escalas temporais.

\begin{ex}
\label{ex:15}
Ao considerar-se a escala temporal $\mathbb{T} = \mathbb{Z}$
e dadas as fun\c{c}\~{o}es $f,g: \mathbb{Z} \rightarrow \mathbb{Z}$
definidas por
\[
f\left(  t\right)  =t^{3}\text{\quad e\quad}g\left(  t\right)  =2t
\]
repara-se que
\[
\left(  f\circ g\right)  ^{\Delta}=\left(  8t^{3}\right)  ^{\Delta}=8\left(
3t^{2}+3t+1\right)  =24t^{2}+24t+8
\]
e que
\[
f^{\Delta}\left(  g\left(  t\right)  \right)  \cdot g^{\Delta}\left(
t\right)  =\left(  12t^{2}+6t+1\right)  2=24t^{2}+12t+2\, .
\]
Logo,
\[
\left(  f\circ g\right)  ^{\Delta}\neq f^{\Delta}\left(  g\left(  t\right)
\right)  g^{\Delta}\left(  t\right)  \, .
\]
\end{ex}

Uma das poss\'{\i}veis maneiras de determinar a derivada da fun\c{c}\~{a}o composta numa
escala temporal \'{e} dada pelo pr\'{o}ximo resultado.

\begin{theorem}[\cite{Bohner2001}]
Seja $g: \mathbb{R} \rightarrow \mathbb{R}$
uma fun\c{c}\~{a}o cont\'{\i}nua tal que $g:\mathbb{T}\rightarrow \mathbb{R}$
\'{e} delta-diferenci\'{a}vel em $\mathbb{T}^{\kappa}$ e seja
$f: \mathbb{R} \rightarrow \mathbb{R}$
uma fun\c{c}\~{a}o com derivada cont\'{\i}nua. Ent\~{a}o, existe um n\'{u}mero real $c\in\left[
t,\sigma\left(  t\right)  \right]  $ tal que
\begin{equation}
\label{eq:2}
\left(  f\circ g\right)  ^{\Delta}\left(  t\right)  =f^{\prime}\left(
g\left(  c\right)  \right)  g^{\Delta}\left(  t\right) \, .
\end{equation}
\end{theorem}

\begin{ex}
\label{ex:exAnt:dfc}
Considere-se as fun\c{c}\~{o}es $f: \mathbb{R} \rightarrow \mathbb{R}$
e $g:\mathbb{T}\rightarrow \mathbb{R}$ definidas por
\[
f\left(  t\right)  =t^{3}\text{\quad e\quad}g\left(  t\right)  =2t \, .
\]
Observe-se que
\begin{equation*}
  \left(  f\circ g\right)  ^{\Delta}\left(  0\right)  =f^{\prime}\left(
g\left(  c\right)  \right)  g^{\Delta}\left(0\right)
  \Leftrightarrow c=\frac{1}{\sqrt{3}}\vee c=-\frac{1}{\sqrt{3}}
\end{equation*}
e, portanto, a f\'{o}rmula \eqref{eq:2} \'{e} v\'{a}lida em $t = 0$
para $c=\frac{1}{\sqrt{3}}\in\left[
0,\sigma\left(  0\right)  \right]  =\left[  0,1\right]$.
\end{ex}

Outra formula\c{c}\~{a}o poss\'{\i}vel para a derivada da fun\c{c}\~{a}o composta \'{e} apresentada no
teorema seguinte.

\begin{theorem}[\cite{Bohner2001}]
Seja $f: \mathbb{R} \rightarrow \mathbb{R}$
uma fun\c{c}\~{a}o diferenci\'{a}vel com continuidade
e $g:\mathbb{T}\rightarrow \mathbb{R}$
uma fun\c{c}\~{a}o com derivada delta. Ent\~{a}o,
\[
\left(  f\circ g\right)  ^{\Delta}\left(  t\right)  =g^{\Delta}\left(
t\right)  \int_{0}^{1}f^{\prime}\left(  g\left(  t\right)  +\xi\mu\left(
t\right)  g^{\Delta}\left(  t\right)  \right)  d\xi\, .
\]
\end{theorem}

\begin{ex}
Considerando as fun\c{c}\~{o}es $f$ e $g$ do Exemplo~\ref{ex:exAnt:dfc}, verifica-se que
\begin{align*}
\left(  f\circ g\right)  ^{\Delta}\left(  t\right)   &= 2\int_{0}^{1}3\left(
2t+2\xi\right)  ^{2}d\xi\\
&= \left[  \left(  2t+2\xi\right)  ^{3}\right]  _{0}^{1}\\
&= \left(  2t+2\right)  ^{3}-8t^{3}\\
&= 24t^{2}+24t+8\, .
\end{align*}
\end{ex}

\begin{theorem}[\cite{Bohner2001}]
\label{teo36}
Seja $g:\mathbb{T}\rightarrow \mathbb{R}$
uma fun\c{c}\~{a}o estritamente crescente tal que $\tilde{\mathbb{T}}=g\left(
\mathbb{T}\right)$ \'{e} uma escala temporal e seja
$f: \tilde{\mathbb{T}}\rightarrow \mathbb{R}$ uma fun\c{c}\~{a}o.
Se $g^{\Delta}\left(  t\right)$ e $f^{\tilde{\Delta}}\left(g\left(t\right)\right)$
existem para $t\in\mathbb{T}^{\kappa}$, ent\~{a}o
\[
\left(  f\circ g\right)  ^{\Delta}=\left(  f^{\tilde{\Delta}}\circ g\right)
g^{\Delta}
\]
onde $f^{\tilde{\Delta}}$ representa a delta derivada de $f$ em
$\tilde{\mathbb{T}}$.
\end{theorem}

\begin{ex}
Considerem-se as fun\c{c}\~{o}es $f$ e $g$ do Exemplo~\ref{ex:15}. Assim $\tilde{\mathbb{T}}
=g\left(\mathbb{Z}\right)  =2 \mathbb{Z}$. Como
\begin{align*}
f^{\tilde{\Delta}}\left(  t\right)   &  =\frac{f\left(  t+2\right)  -f\left(
t\right)  }{2}\\
&  =\frac{\left(  t+2\right)  ^{3}-t^{3}}{2}\\
&  =3t^{2}+6t+4
\end{align*}
e
\[
f^{\tilde{\Delta}}\circ g\left(  t\right)  =12t^{2}+12t+4
\]
logo
\[
\left(  f^{\tilde{\Delta}}\circ g\right)  g^{\Delta}\left(  t\right)  =\left(
12t^{2}+12t+4\right)  2=24t^{2}+24t+8\, .
\]
\end{ex}

Como consequ\^{e}ncia do Teorema~\ref{teo36} obt\'{e}m-se o seguinte resultado
para a derivada da fun\c{c}\~{a}o inversa.

\begin{theorem}[\cite{Bohner2001}]
Seja $g:\mathbb{T}\rightarrow \mathbb{R}$
uma fun\c{c}\~{a}o estritamente crescente tal que
$\tilde{\mathbb{T}} =g\left(\mathbb{T}\right)$ \'{e} uma escala temporal. Ent\~{a}o
\[
\frac{1}{g^{\Delta}}=\left(  g^{-1}\right)  ^{\tilde{\Delta}}\circ g
\]
nos pontos $t\in\mathbb{T}^\kappa$ onde $g^{\Delta}\left(  t\right)  \neq0$.
\end{theorem}


\newpage

\section{Conclus\~{a}o e notas finais}
\label{sec:4}

O formalismo das escalas temporais permite desenvolver
uma teoria de c\'{a}lculo diferencial generalizado,
que tem como casos particulares o c\'{a}lculo diferencial em $\mathbb{R}$,
o c\'{a}lculo das diferen\c{c}as finitas e o c\'{a}lculo-$q$ (qu\^{a}ntico).
Mesmo uma escala temporal como o conjunto de Cantor
est\'{a} inclu\'{\i}da. Os objectivos primordiais da an\'{a}lise em escalas temporais
s\~{a}o \emph{unifica\c{c}\~{a}o} e \emph{generaliza\c{c}\~{a}o}. Como resultado
pode-se evitar a apresenta\c{c}\~{a}o em paralelo de resultados discretos e cont\'{\i}nuos
(por vezes ma\c{c}ador, \`{a}s vezes dif\'{\i}cil) e a realiza\c{c}\~{a}o de ``demonstra\c{c}\~{o}es''
por analogia ou apenas por omiss\~{a}o, que t\~{a}o frequentemente ocorrem
nas apresenta\c{c}\~{o}es de vers\~{o}es discretas de temas da an\'{a}lise \cite{Hilger1990}.

Uma \'{a}rea de investiga\c{c}\~{a}o muito activa
consiste em estudar sistemas de controlo em escalas temporais
(\emph{vide}, \textrm{e.g.}, \cite{Bartosiewicz2007,B:Paw,Bartosiewicz2008,D2,withEwaP:avoidance,Gosia:Delfim}).
Este \'{e} o assunto da unidade curricular \emph{Teoria do Controlo em Time Scales}
da \'{a}rea de especializa\c{c}\~{a}o em \emph{Optimiza\c{c}\~{a}o, Sistemas e Controlo}
do \emph{Programa Doutoral em Matem\'{a}tica e Aplica\c{c}\~{o}es} (PDMA)
entre os Departamentos de Matem\'{a}tica da Universidade de Aveiro (DMAT-UA)
e da Universidade do Minho (DMAT-UM). De modo muito sucinto,
a ideia central da teoria dos sistemas e controlo em \emph{time scales}
consiste em usar o conceito de derivada delta (ou nabla) para unificar
os sistemas de controlo tanto em tempo cont\'{\i}nuo como em tempo discreto:

\bigskip

\begin{minipage}[h]{\MiniPageLeft}
\begin{center}
   {\it Tempo  cont\'{\i}nuo}\\
   $t\in \mathbb{R}$
\end{center}
\vspace{-5mm}
\[\rnode{CCSR}{\psframebox[linecolor=black]{\dot{x}(t)=f_1(x(t),u(t))}}\]
Se $\mathbb{T}=\mathbb{R}$, ent\~{a}o
\[\dot{x}(t)=x^\Delta(t).\]
Logo,
\[\rnode{CCST}{\psframebox[linecolor=black]{x^\Delta (t)=f_1(x(t),u(t))}}\]
\nccurve[linecolor=black,linestyle=dotted,angleA=270,angleB=90]{->}{CCSR}{CCST}

\end{minipage}
\begin{minipage}[h]{\MiniPageRight}
\begin{center}
    {\it Tempo discreto}\\
    $t\in \mathbb{Z}$
\end{center}
\vspace{-5mm}
\[\rnode{DCSZ}{\psframebox[linecolor=black]{x(t+1)=f_2(x(t),u(t))}}\]
Se $\mathbb{T}=\mathbb{Z}$, ent\~{a}o
\[x^\Delta (t)=x(t+1)-x(t).\]
Logo,
\[\rnode{DCST}{\psframebox[linecolor=black]{x^\Delta (t)=f_2(x(t),u(t))-x(t)}}\]
\nccurve[linecolor=black,linestyle=dotted,angleA=270,angleB=90]{->}{DCSZ}{DCST}
\end{minipage}

Estudam-se ent\~{a}o sistemas de controlo numa escala temporal $\mathbb{T}$, \textrm{i.e.},
\begin{equation*}
x^\Delta (t)=f(x(t),u(t)),\quad t\in \mathbb{T}^\kappa \, .
\end{equation*}
Ao considerar-se o c\'{a}lculo em escalas temporais,
os sistemas de controlo em tempo cont\'{\i}nuo e em tempo discreto
s\~{a}o fundidos numa \'{u}nica teoria mais geral.

Neste pequeno artigo apresentam-se as no\c{c}\~{o}es elementares da teoria das escalas temporais,
que servem de ponto de partida para um estudo mais aprofundado,
\textrm{e.g.}, o estudo do c\'{a}lculo integral
em escalas temporais, o estudo das equa\c{c}\~{o}es delta ou nabla-diferenci\'{a}veis
e a ``Teoria do Controlo em Time Scales''.
O \emph{package} \emph{TimeScales}, desenvolvido para o \textsf{Mathematica},
n\~{a}o s\'{o} executa c\'{a}lculos b\'{a}sicos como evidencia as diferen\c{c}as de resultados
para escalas temporais distintas. Est\'{a} dispon\'{\i}vel em
\url{http://www.esce.ipvc.pt/docentes/srodrigues/TimeScales.m}.

Em anos recentes o c\'{a}lculo em escalas temporais
tem recebido consider\'{a}vel aten\c{c}\~{a}o internacional
(\textrm{cf.}, \textrm{e.g.}, \url{http://web.mst.edu/~bohner/tslist.html}).
V\'{a}rios semin\'{a}rios e sess\~{o}es especiais sobre a tem\'{a}tica das
escalas temporais s\~{a}o organizados.
Em Maio pr\'{o}ximo decorre em Dresden, na Alemanha, uma sess\~{a}o convidada
organizada por Martin Bohner, Stefan Hilger e Agacik Zafer,
no \^{a}mbito da \emph{8th AIMS Conference on Dynamical Systems, Differential Equations and Applications},
Dresden, Germany, May 25-28, 2010 (\url{http://web.mst.edu/~bohner/dd2010.html}).
Ao leitor interessado nestas tem\'{a}ticas, e como ponto de partida,
recomenda-se os livros \cite{Bohner2001,Bohner2003,Lak:book}.


\section*{Agradecimentos}

Os dois primeiros autores s\~{a}o alunos do
Programa Doutoral em Matem\'{a}tica e Aplica\c{c}\~{o}es (PDMA)
dos Departamentos de Matem\'{a}tica das Universidades de Aveiro (DMAT-UA)
e Minho (DMAT-UM), com suporte financeiro
da Funda\c{c}\~{a}o para a Ci\^{e}ncia e Tecnologia (FCT):
bolsas SFRH/BD/33634/2009 (Artur Cruz) e SFRH/BD/33384/2008 (Helena Rodrigues).
Os autores est\~{a}o gratos a Pedro A. F. Cruz pela ajuda na constru\c{c}\~{a}o
do \emph{package} em \textsf{Mathematica} \emph{TimeScales},
a Nat\'{a}lia Martins por todo o acompanhamento
e pela leitura cuidada de uma vers\~{a}o preliminar deste trabalho
e a um revisor an\'{o}nimo pelas in\'{u}meras sugest\~{o}es de melhoramento do texto.



\end{document}